\RequirePackage{fix-cm}
\documentclass[11pt,a4paper,leqno]{article}
\usepackage{graphicx}
\usepackage{mathptmx}   

\usepackage[latin1]{inputenc}  
\usepackage{amsfonts}       
\usepackage{amsmath}
\usepackage{amssymb}
\usepackage{multirow}
\usepackage{graphicx}  
\usepackage{enumerate}
\usepackage{bibentry}
\usepackage{enumerate}
\usepackage{color}
\usepackage{amsmath}

\newcommand{\I}{\mathbb{I}}
\newcommand{\Ih}{\mathbb{I}_h}
\newcommand{\eqdist}{\stackrel{d}{=}}

\usepackage{amsthm}
\newtheorem{thm}{Theorem}

\theoremstyle{definition}
\newtheorem{ex}{Example}
\newtheorem{definition}{Definition}
\newtheorem{remark}{Remark}

\newtheorem{lemma}{Lemma}

\begin{document}

\title{A comparison method for heavy-tailed random variables}

\author{Jaakko Lehtomaa\footnote{e-mail address: jaakko.lehtomaa@helsinki.fi}}

\maketitle

\begin{abstract}
We investigate a way of comparing and classifying tails of random variables. Our approach extends the notion of classical indices, such as exponential and moment indices, which are widely used measuring heaviness of tail functions. A non-parametric risk measure applicable for all heavy-tailed random variables is obtained as a concave function that represents the decay speed of tail function. Many key properties of the distribution of a random variable are encoded into this function, which enables a new way to estimate tails. The latter half of the paper is devoted to numerous examples illustrating properties of the results developed during the first half. 
\end{abstract}


\section{Introduction}\label{intro}

Suppose $(\Omega, \mathcal{F},P)$ is a probability space where subsequent random variables are defined. For a random variable $X$, with distribution function $F_X(x)=P(X\leq x)$ and tail function $\overline{F}_X=1-F_X$, we define the  \emph{hazard function} by $R_X=-\log \overline{F}_X$. All random variables are assumed essentially unbounded from above, that is, $P(X>a)>0$ for all $a>0$. 

Let $h \colon [0,\infty) \to [0,\infty)$ be an increasing function with the property $h(x) \to \infty$, as $x \to \infty$. Such a function $h$ is called a \emph{scale function}. The quantity
\begin{equation*}
\liminf_{x \to \infty} \frac{R_X(x)}{h(x)}
\end{equation*}
reveals information about the integrability properties of $X$ or, more precisely, about the integrability of the transformation $h(X)$. Choosing $h(x)=x$ or $h(x)=\log x$ (although this is not a scale function) one obtains the classical \emph{exponential index} $\mathcal{E}(X)$ or \emph{moment index} $\mathbb{I}(X)$ respectively. Moreover, with these choices interpretations 
\begin{equation}\label{e}
\mathcal{E}(X)=\liminf_{x \to \infty} \frac{R_X(x)}{x}=\sup\{s\geq 0: E(e^{sX})<\infty\}
\end{equation}
and
\begin{equation}\label{m}
\mathbb{I}(X)=\liminf_{x \to \infty} \frac{R_X(x)}{\log x}=\sup\{s\geq 0: E((X^+)^s)<\infty\},
\end{equation}
where $x^+=\max(0,x)$, are valid as can be seen from e.g. \cite{daley2001moment} and \cite{daley2006moment}. 

We will mainly study non-negative random variables. However, most of the properties can be transferred to the unrestricted case simply by considering variable $X^+$ instead of $X$, since these two variables have the same right tail. To study left tails, one can replace $X$ by $-X$.

\subsection{The main problem and proposed solution}

The difficulty with indices defined in formulas \eqref{e} and \eqref{m} is that neither of them can be used to compare tails of random variables $X$ and $Y$ if their indices share the same value $0$ or $\infty$. In this situation the scale $h$ does not represent correctly the scale of the hazard functions $R_X$ and $R_Y$. This raises two questions:
\begin{enumerate}[\bf{Q}1.]
\item \label{q1} Given two general random variables, how can their tails be compared?
\item \label{q2} How could one measure the heaviness of a general heavy-tailed random variable?
\end{enumerate} 

It seems that questions \ref{q1} and \ref{q2} have not been studied extensively in the past. However, these kinds of questions have recently attracted attention among practitioners of risk management. In \cite{halliwell2013classifying} one can find an applied approach with related discussion to the tail comparison problem. We will provide a completely different solution that is applicable to a wider class of probability distributions. 

To answer the question \ref{q1}, we propose a direct comparison between the associated hazard functions via quantity
\begin{equation}\label{ratio}
\liminf_{x \to \infty} \frac{R_X(x)}{R_Y(x)}.
\end{equation}
If the quantity of formula \eqref{ratio} equals $a\in(0,\infty)$, we may deduce that for any small $\epsilon>0$ there exists a number $x_\epsilon$ such that for all $x>x_\epsilon$ the inequality
\begin{equation*}
\overline{F}_X(x)\leq \overline{F}_Y(x)^{a-\epsilon}
\end{equation*}
holds. This enables comparison between $X$ and $Y$ even when indices \eqref{e} and \eqref{m} fail to characterise the proper decay speed.

In addition to the direct comparison of type \eqref{ratio}, it will be shown that any risk function of a heavy-tailed random variable can, in a sense, be replaced by a suitable concave scale function that adequately represents its asymptotic scale. This answers the question \ref{q2}: Heaviness is measured by the asymptotic growth speed of this deterministic function. Using a concave function is beneficial because it has, in many cases, much simpler  representation than the original hazard function. 

\subsection{Structure of the paper}

The rest of the paper is arranged as follows. In chapter \ref{motivation} necessary background information is given with preliminary results. In chapter \ref{results}, we develop the main properties that are used in the applications chapter \ref{examples}. Lastly, technical constructions omitted in section \ref{results} are given in appendix chapter \ref{mappendix}.

\section{Motivation}\label{motivation}
A random variable with a positive exponential index is called \emph{light-tailed}. For such variables one can always deduce that the speed at which the tail function $\overline{F}$ decreases is at least exponential, that is, 
\begin{equation}\label{exp}
\overline{F}(x)\leq e^{-ax}
\end{equation}
for some $a>0$ and all  $x$ large enough.
If a random variable is not light-tailed, it is called \emph{heavy-tailed}. In the case of a positive and finite moment index  polynomial bound
\begin{equation*}
\overline{F}(x)\leq x^{-b}
\end{equation*}
for some $b>0$ and all large enough $x$ may be obtained, whereas inequality of the type \eqref{exp} is not possible. We aim to provide a bound suitable for all heavy-tailed random variables in the form 
\begin{equation}\label{general}
\overline{F}(x)\leq e^{-h(x)}
\end{equation}
for all $x$ large enough, where the function $h$ is an accurate representation of the true decay speed of the tail $\overline{F}$. In example \ref{esim2} we will see how this can be achieved in the case of Weibull or log-normal type distributions. In order to find a suitable function $h$ we introduce the following definition.

\begin{definition}\label{ekadef} Suppose $X$ is a random variable and $h$ is a scale function. Then 
\begin{equation*}
\Ih(X)=\liminf_{x \to \infty} \frac{R_{X}(x)}{h(x)}
\end{equation*}
is called the $h$-order of X. If $\Ih(x)=1$, the function $h$ is called a \emph{natural scale (function)} of $X$. Hereafter, $h_X$ denotes a natural scale of $X$.
\end{definition}

\begin{remark} In definition \ref{ekadef}, the concept of natural scale does not uniquely define any function $h$. Instead, there are many different choices. Trivial candidate is always $h=R_X=-\log \overline{F}_X$. This, however, may turn out to be a cumbersome choice. 
\end{remark}

Suppose $h_X$ is a natural scale of a random variable $X$, then for every $\epsilon>0$ inequality similar to \eqref{general} holds:
\begin{equation}\label{esittely}
\overline{F}(x)\leq e^{-(1-\epsilon)h_X(x)}
\end{equation}
for all $x$ large enough.  We will see in Theorem \ref{theorem3} that, for a heavy-tailed random variable, a function $h_X$ can always be chosen so that:
\begin{enumerate}[I]
\item $h_X$ is concave\label{omin1}
\item $h_X(0)=0$\label{omin2}
\item $h_X$ is essentially the best choice in \eqref{esittely} \label{omin3}.
\end{enumerate}
Properties \ref{omin1} and \ref{omin2} can make the function $h_X$ smoother than the original risk function $R$ itself. However, $h_X$ conveys useful information about the asymptotic behaviour of the tail $\overline{F}$. This becomes apparent when studying expectations. We will see that there exist numbers $a,b \in (0,\infty)$, where $a<b$, such that
\begin{equation}\label{111}
E(e^{ah_X(X)})<\infty,
\end{equation}
but
\begin{equation}\label{222}
E(e^{bh_X(X)})=\infty.
\end{equation}
This means that the function $h_X$ regularises the random variable $X$ so that the expectation \eqref{111} is finite, but sparingly enough for expectation \eqref{222} to be divergent. Interpretation of this is that the deterministic function $h_X$ captures the scale of the random variable $X$ and thus  measures how risky the variable is. Precise information about the magnitude of the possibility of very large realisations is of crucial importance in many fields. For example, in insurance and finance large losses are possible, say, in catastrophe insurance or in derivatives trading.

Theorems for general moment properties of heavy-tailed random variables can be found from chapter 2 of  \cite{foss2011introduction}. We conclude the chapter by recalling one of these theorems.

\begin{thm}[Theorem $2.9$ of \cite{foss2011introduction}] \label{FKZ} Let $X\geq 0$ be a heavy-tailed random variable. Suppose $g$ is a real valued function for which $g(x) \to \infty$, as $x \to \infty$. 

Then, there exists a monotone concave function $h$ satisfying properties:
\begin{enumerate}
\item $h(x)=o(x)$, as $x \to \infty$ \label{l1}
\item $E(e^{h(X)})<\infty$ \label{l2}
\item $E(e^{h(X)+g(X)})=\infty$.\label{l3}
\end{enumerate}
\begin{proof}
See \cite{foss2011introduction} pp. 12-13.
\end{proof}
\end{thm}

From now on, we will omit the lower index indicating the random variable from hazard and tail functions whenever the variable in question is clear from context. In addition, sequences will be denoted in short as $(x_n):=(x_n)_{n=1}^\infty$, where $:=$ denotes equality by definition.

\section{Main results}\label{results}

\subsection{Existence of a suitable natural scale}

One of the main results of the paper is Theorem \ref{theorem3}, where the existence of a desired function $h$ satisfying requirements \ref{omin1}-\ref{omin3} of section \ref{motivation} is shown. Its proof requires the next theorem, which translates a relation similar to \eqref{e} or \eqref{m} into a more general environment.

\begin{thm}\label{theorem1} Suppose $X$ is a random variable. Assume further that $h$ is continuous and $h(x)\to \infty$ as $x \to \infty$.  Then 
\begin{equation}\label{rep1}
\liminf_{x \to \infty} \frac{R(x)}{h(x)}=\sup \{ s\geq 0: E(e^{sh(X)})<\infty\}.
\end{equation}
\begin{proof} We divide the proof in two parts.
\begin{enumerate} \item \label{ekaosa}
Suppose first that the function $h$ is strictly increasing. Applying \eqref{e} to random variable $h(X)$ yields
\begin{equation*}
\sup \{ s\geq 0: E(e^{sh(X)})<\infty\} =-\limsup_{x \to \infty} \frac{\log P(h(X)>x)}{x}.
\end{equation*}
Since $h$ is invertible, we obtain
\begin{eqnarray*}
\limsup_{x \to \infty} \frac{\log P(h(X)> x)}{x}&=&\limsup_{x \to \infty} \frac{\log P(X> h^{-1}(x))}{x} \\
&=& \limsup_{z \to \infty} \frac{\log P(X> h^{-1}(h(z)))}{h(z)} \\
&=& \limsup_{z \to \infty} \frac{\log P(X> z)}{h(z)}.
\end{eqnarray*}
This ends the proof of part \ref{ekaosa}.

\item \label{tokaosa} Suppose then that the function $h$ is increasing, but not necessarily strictly increasing. Let $\eta>0$. We may choose a strictly increasing continuous function $h_\eta$ such that for all $x\geq 0$ 
\begin{equation}\label{tuplayht}
h(x)\leq h_\eta(x)\leq h(x)+\eta
\end{equation}
holds. See appendix \ref{eappendix} below for the actual construction of function $h_\eta$. By part \ref{ekaosa} the result \eqref{rep1} holds for the function $h_\eta$. Using \eqref{tuplayht} it is easy to see that 
$$\liminf_{x \to \infty} \frac{R(x)}{h_\eta(x)}=\liminf_{x \to \infty} \frac{R(x)}{h(x)} $$
and 
$$\sup \{ s\geq 0: E(e^{sh_\eta(X)})<\infty\}=\sup \{ s\geq 0: E(e^{sh(X)})<\infty\}, $$
which ends the proof.
\end{enumerate}
\end{proof}
\end{thm} 

We are now in a position to show that a natural scale, defined in definition \ref{ekadef}, can always be chosen in  the following way.

\begin{thm}\label{theorem3} Suppose $X$ is a heavy-tailed random variable. Then there exists a concave function $h$ for which $h(0)=0$ such that 
\begin{equation}\label{hyvakaava}
\liminf_{x \to \infty} \frac{R(x)}{h(x)}=1.
\end{equation}
Equivalently, there exist numbers $a,b \in (0,\infty)$ and a concave function $h^*$ for which $h^*(0)=0$ such that 
\begin{equation}\label{kas}
E(e^{ah^*(X)})<\infty
\end{equation}
and
\begin{equation}\label{kas2}
E(e^{bh^*(X)})=\infty.
\end{equation}
\begin{proof} Equivalence of the above assertions is immediate. If \eqref{hyvakaava} holds, we may choose $h^*=h$, $a=1/2$ and $b=3/2$ in \eqref{kas} and \eqref{kas2}. Result is implied by Theorem \ref{theorem1}. For the other direction, Theorem \ref{theorem1} tells us that $\alpha:=\liminf_{x \to \infty} R(x)/h^*(x)\in [a,b]$. Setting $h=\alpha h^*$ gives the required function.

We will thus concentrate on proving formula \eqref{hyvakaava}. To see this, let $g$ be a non-negative continuous function for which $g(x)=o(R(x))$ holds, as $x \to \infty$. For explicit construction of such a function see appendix \ref{tappendix} below. Now, there exists a function $\hat{h}$ satisfying conditions \ref{l1}-\ref{l3} of Theorem \ref{FKZ}. 

Note that the function $\hat{h}$ cannot be bounded from above. If $\hat{h}$ was bounded by a positive constant $M$, we would get 
\begin{equation}\label{aarvio}
E(e^{\hat{h}(X)+g(X)})\leq e^{M}E(e^{g(X)}).
\end{equation}
However, because $\liminf_{x \to \infty} R(x)/g(x)=\infty$, Theorem \ref{theorem1} shows that especially $E(e^{g(X)})<\infty$ holds. Formula \eqref{aarvio} combined with requirement \ref{l3} of Theorem \ref{FKZ} would now imply a contradiction.

Using Theorem \ref{theorem1} and properties \ref{l2} and \ref{l3} of Theorem \ref{FKZ} we get
\begin{equation}\label{yk1}
\liminf_{x \to \infty} \frac{R(x)}{\hat{h}(x)}\geq 1
\end{equation}
and
\begin{equation}\label{yk2}
\liminf_{x \to \infty} \frac{R(x)}{\hat{h}(x)+g(x)}\leq 1.
\end{equation}
If the relation $\hat{h}(x)=o(R(x))$ was valid, it would imply $\hat{h}(x)+g(x)=o(R(x))$, which in turn would imply $\lim_{x \to \infty} R(x)/(\hat{h}(x)+g(x))=\infty$. This, however, is impossible by the statement of formula \eqref{yk2}. Hence, relation $\hat{h}(x)=o(R(x))$ cannot hold, which implies $\liminf_{x \to \infty} R(x)/\hat{h}(x)<\infty$.

Finally, denoting 
$$\beta:=\liminf_{x \to \infty} R(x)/\hat{h}(x)\in [1,\infty)$$
and setting 
$$h(x):=\beta (\hat{h}(x)-\hat{h}(0))$$
gives the desired function.
\end{proof}
\end{thm}

\begin{remark}\label{remark1} If $h$ is a concave function with $h(0)=0$, the subadditivity requirement
$$h(a+b)\leq h(a)+h(b) $$
holds for all $a,b>0$. In addition, it is easy to check that for any $a \in (0,1)$ relation $\limsup_{x \to \infty} h(ax)/h(x)< \infty$ holds. This means, in particular, that the function $h$ belongs to dominatedly varying class $\mathcal{D}$. See appendix  \ref{rappendix} below for details. Moreover, if $h$ is also a scale function, it satisfies the asymptotic relation $h(x)\to \infty$, as $x \to \infty$. This implies that $h$ is continuous and strictly increasing.
\end{remark}

Theorem \ref{theorem3} gives a way to classify random variables purely by the thickness of their tails. This thinking is different from many other classifications of heavy-tailed random variables where an analytic property, not explicitly related to the tail decay speed, is required. 

\begin{remark}
Theorem \ref{theorem3} shows how to find a natural scale for a heavy-tailed random variable. Namely, if a concave function satisfying \eqref{kas} and \eqref{kas2} is found, it is a natural scale up to a positive multiplicative constant. Good initial guess for finding a suitable function $h$ is $R=-\log \overline{F}$ or a suitable dominating component of $R$. In example \ref{esim4} we will illustrate properties of this choice.
\end{remark}

\subsection{Properties of natural scales}
Properties of indices \eqref{e} and \eqref{m} are different. For example, if $X$ and $Y$ are positive and independent, the equality $\mathbb{I}(XY)=\min(\mathbb{I}(X),\mathbb{I}(Y))$ is always valid whereas $\mathcal{E}(XY)=\min(\mathcal{E}(X),\mathcal{E}(Y))$ is not (e.g. take constant $Y$). 

We provide necessary conditions ensuring that the $h$-order of the sum and product of independent variables is the minimum of the associated $h$-orders. These properties allow one to make simple and fast estimates even if the exact computation in not feasible. The next theorem gives sufficient conditions for simple computational rules to hold. The aim is to establish results that can be tested with natural scales of random variables.

\begin{thm}\label{laskutheorem} Suppose $X$ and $Y$ are positive, independent and essentially unbounded random variables. Assume that $h$ is a continuous function and $h(x)\to \infty$, as $x\to \infty$. 

Then, the following implications hold:
\begin{equation}\label{summa}
\forall a,b>0: \, h(a+b)\leq h(a)+h(b) \Longrightarrow  \Ih(X+Y)=\min(\Ih(X),\Ih(Y))
\end{equation}
and 
\begin{equation}\label{tulo}
\forall a,b>0: \, h(ab)\leq h(a)+h(b) \Longrightarrow \Ih(XY)=\min(\Ih(X),\Ih(Y)).
\end{equation}

\begin{proof} We use the representation \eqref{rep1} of Theorem \ref{theorem1} with the facts 
\begin{equation*}
E(e^{sh(X+Y)}) \leq E(e^{sh(X)})E(e^{sh(Y)})
\end{equation*}
and  
\begin{equation*}
E(e^{sh(XY)}) \leq E(e^{sh(X)})E(e^{sh(Y)})
\end{equation*}
from formulas \eqref{summa} and \eqref{tulo} to see that $\Ih(X+Y)\geq \min(\Ih(X),\Ih(Y))$ and $\Ih(XY)\geq \min(\Ih(X),\Ih(Y))$ are implied. For converse inequalities one may again use \eqref{rep1} combined with positivity to see that 
\begin{equation}\label{ee1}
E(e^{sh(X+Y)})\geq E(e^{sh(X)})
\end{equation}
and
\begin{equation}\label{ee2}
E(e^{sh(XY)})\geq E(e^{sh(XY)}\mathbf{1}(Y\geq 1))\geq E(e^{sh(X)})P(Y\geq 1).
\end{equation}
Formulas \eqref{ee1} and \eqref{ee2} imply $\Ih(X+Y)\leq \Ih(X)$ and $\Ih(XY)\leq \Ih(X)$. The remaining case is clear by the symmetry of random variables $X$ and $Y$.
\end{proof}
\end{thm}

\begin{remark}\label{cor1} If $X$ and $Y$ are heavy-tailed random variables, Theorem \ref{theorem3} ensures that natural scales $h_X$ and $h_Y$ can be chosen concave. Thus, by Theorem \ref{laskutheorem} and remark \ref{remark1} we may use the calculation rule 
$$\Ih(X+Y)=\min(\Ih(X),\Ih(Y)),$$
where $h=h_X$ or $h=h_Y$.
\end{remark}

Purpose of the Theorem \ref{laskutheorem} is to give conditions that enable simple calculations. This is why random variables are assumed independent. However, the following result confirms that in certain cases we may infer scales even without independence.

\begin{thm}\label{vertailuthm} Suppose $X$ and $Y$ are positive heavy-tailed random variables. Let $h_X$ and $h_Y$ be concave natural scales of $X$ and $Y$ with $h_Y(0)=0$ (obtained e.g. from Theorem \ref{theorem3}). Assume further that 
\begin{equation}\label{roletus}
\lim_{x \to \infty} \frac{h_X(x)}{h_Y(x)}=\infty.
\end{equation}

Then there exists $c \in (0,\infty)$ such that
\begin{equation}\label{summaskaala}
h_{X+Y}(x)=c h_Y(x)
\end{equation}
is a natural scale of $X+Y$.
\begin{proof} Because of \eqref{roletus} is it clear that 
\begin{equation*}
\liminf_{x \to \infty} \frac{R_X(x)}{h_Y(x)}\geq \left( \liminf_{x \to \infty} \frac{R_X(x)}{h_X(x)} \right) \left(\liminf_{x \to \infty}\frac{h_X(x)}{h_Y(x)} \right)=\infty,
\end{equation*}
which implies, using Theorem \ref{theorem1}, that $E(e^{sh_Y(X)})<\infty$ for all $s>0$. Now, since $X$ and $Y$ are assumed positive, 
\begin{equation*}
E(e^{sh_Y(X+Y)})\geq E(e^{sh_Y(Y)}).
\end{equation*}
This implies $\I_{h_Y}(X+Y)\leq \I_{h_Y}(Y)=1$. On the other hand, because $h_Y$ is by remark \ref{remark1} subadditive, we get
\begin{eqnarray*}
& &E(e^{sh_Y(X+Y)}) \\
&\leq & E(e^{sh_Y(X)+sh_{Y}(Y)}) \\
& = & E(e^{sh_Y(X)+sh_{Y}(Y)}\mathbf{1}(X\geq Y))+E(e^{sh_Y(X)+sh_{Y}(Y)}\mathbf{1}(X< Y)) \\
&\leq & E(e^{2sh_Y(X)})+E(e^{2sh_Y(Y)}).
\end{eqnarray*}
This implies $\I_{h_y}(X+Y)\geq (1/2)\I_{h_Y}(Y)=1/2$. Definition of natural scale yields now \eqref{summaskaala}.
\end{proof}
\end{thm}

The last theorem allows one to estimate tails of transformations of IID (independent and identically distributed) variables using the tail of a single variable. This estimate is useful in the study of products, see examples \ref{2prod}, \ref{tasymp} and \ref{aesim} below. Before this result, we need a lemma that expands a central property of indices \eqref{e} and \eqref{m}. 
\begin{lemma}\label{maxlemma} Suppose $X$ and $Y$ are random variables and $h$ is a scale function. Then 
\begin{equation*}
\liminf_{x \to \infty} \frac{R_{\max(X,Y)}(x)}{h(x)}=\min \left( \liminf_{x \to \infty} \frac{R_{X}(x)}{h(x)},\liminf_{x \to \infty} \frac{R_{Y}(x)}{h(x)} \right).
\end{equation*}
\begin{proof} Observe that 
\begin{eqnarray*}
\limsup_{x \to \infty} \frac{\log \overline{F}_{\max(X,Y)}(x)}{h(x)}& \leq & \limsup_{x \to \infty} \frac{\log (\overline{F}_{X}(x)+\overline{F}_{Y}(x))}{h(x)} \\
&=&\max\left( \limsup_{x \to \infty} \frac{\log \overline{F}_{X}(x)}{h(x)},\limsup_{x \to \infty} \frac{\log \overline{F}_{Y}(x)}{h(x)}\right),
\end{eqnarray*}
where the last equality follows from e.g. \cite{dembo2009large} Lemma 1.2.15.
On the other hand,
$$\limsup_{x \to \infty} \frac{\log \overline{F}_{\max(X,Y)}(x)}{h(x)}\geq \limsup_{x \to \infty} \frac{\log \overline{F}_{Z}(x)}{h(x)},$$
where $Z=X$ or $Z=Y$. This proves the claim.
\end{proof}
\end{lemma}

\begin{thm}\label{vtheorem}
Let $n\geq 2$ be fixed and suppose $X,X_1,X_2,\ldots,X_n$ are positive heavy-tailed IID variables with continuous common distribution function $F$. Assume further that $g\colon \mathbb{R}^{n}\to \mathbb{R}$ is a function with properties:
\begin{enumerate}
\item Each component $g_i$, $i\in \{1,2,\ldots,n\}$, of function $g$ is an increasing function and $g_i(x)\to \infty$, as $x \to \infty$.
\item The diagonal function $g_d(x):=g(x,x,\ldots,x)$ has inverse function, denoted by $g_d^{-1}$.
\end{enumerate}

Then there exists a positive constant $c$ such that
\begin{equation}\label{tulokaava1}
h_{g(X_1,X_2,\ldots ,X_n)}(x)=c (R_{X}\circ g_d^{-1})(x)
\end{equation}
is a natural scale of transformation $g(X_1,X_2,\ldots,X_n)$. Moreover, for any $\epsilon>0$ there exists a number $x_\epsilon$ such that for all $x> x_\epsilon$ the following inequality holds:
\begin{equation}\label{arviokaava}
P(g(X_1,X_2,\ldots,X_n)>x)\leq P(X_1>g_d^{-1}(x))^{1-\epsilon}.
\end{equation}
\begin{proof} To prove the claim we will initially set $\hat{h}=(R_{X}\circ g_d^{-1})(x)$. Now
\begin{eqnarray*}
E(e^{s \hat{h}(g(X_1,X_2\ldots, X_n))})&\leq & E(e^{s \hat{h}(g_d(\max(X_1,X_2,\ldots,X_n))}) \\
&=& E(e^{sR_{X}(\max(X_1,X_2,\ldots,X_n))}).
\end{eqnarray*}
Using Lemma \ref{maxlemma} inductively with choice $h=R_{X}$ we see that
\begin{eqnarray*}
\liminf_{x \to \infty} \frac{R_{\max(X_1,X_2,\ldots,X_n)}(x)}{R_{X}(x)}=1.
\end{eqnarray*}
This, together with Theorem \ref{theorem1}, yields 
$$\mathbb{I}_{\hat{h}}(g(X_1,X_2\ldots ,X_n))\geq \I_{R_{X}}(X_1)=1.$$
For the other direction we estimate
\begin{eqnarray}
& & E(e^{s \hat{h}(g(X_1,X_2,\ldots, X_n))}) \nonumber \\
& \geq &E(e^{s \hat{h}(g_d(X_1))}\mathbf{1}(\min(X_2,\ldots,X_n)> X_1)). \label{ekayht1}
\end{eqnarray}
Note that for any positive measurable function $r$ we have
\begin{eqnarray*}
& & E(r(X_1)\mathbf{1}(\min(X_2,\ldots,X_n)> X_1)) \\
&=&\int_{0}^\infty \int_{x_1}^\infty \ldots \int_{x_1}^\infty \, r(x_1) \, \mathrm{d}F(x_n)\mathrm{d}F(x_{n-1})\ldots\mathrm{d}F(x_1) \\
&=&\int_{0}^\infty \, r(x_1) P(X_1>x_1)^{n-1} \, \mathrm{d}F(x_1) \\
&=& E(r(X_1)\overline{F}_{X}(X_1)^{n-1}).
\end{eqnarray*} 
Applying this to formula \eqref{ekayht1} with the choice $r(x)=e^{s \hat{h}(g_d(X_1))}=e^{s R_{X}(X_1)}$ we finally obtain 
\begin{equation*}
E(e^{s R_{X}(X_1)}\overline{F}_{X}(X_1)^{n-1})=E(e^{(s-n+1) R_{X}(X_1)}). 
\end{equation*}
Thus, 
$$\mathbb{I}_{\hat{h}}(g(X_1,X_2,\ldots, X_n))\leq \I_{R_{X}}(X_1)+n-1=n.$$ In conclusion, we deduce that the function $\hat{h}$ differs from the natural scale only by a positive constant factor. Therefore, there exists $c \in [1,n]$ such that \eqref{tulokaava1} holds.

Now, using the definition of natural scale, we get 
$$\liminf_{x \to \infty} \frac{R_{g(X_1,X_2,\ldots,X_n)}(x)}{c (R_{X}\circ g_d^{-1})(x)}=1, $$
from which the claim \eqref{arviokaava} directly follows.
\end{proof}
\end{thm}

\section{Applications and examples}\label{examples}

Suppose we are given two sequences of random variables $(A_i)$ and $(B_i)$. Define 
\begin{equation}\label{sum}
S_n=B_1+B_2+\ldots+B_n
\end{equation}
and 
\begin{equation}\label{disc}
Y_n=B_1+A_1B_2+A_1A_2B_3+\ldots + A_1\ldots A_{n-1} B_n.
\end{equation}
Formula \eqref{sum} could be viewed as the sum of risks from different sources. For example, variables $B_1$, $B_2$ and $B_3$ could represent the aggregate losses of different insurance lines such as casualty, life and catastrophe insurance, where the tail behaviour of each line can be different. Formula \eqref{disc} on the other hand can be viewed as a randomly discounted random cash flow. During the following examples we will see how the $h$-orders of variables $S_n$ and $Y_n$ can be studied using results of the previous chapter. 

We will begin with an example that clarifies why lower limit is used in the definition of $h$-orders.

\begin{ex}[Justification of limes inferior in definition \ref{ekadef}]\label{esim1}
Consider continuous functions $h_1,h_2 \colon [0,\infty) \to [0,\infty)$, where $h_2(x)<h_1(x)$ for all $x\geq 0$. Assume further that the functions $h_1$ and $h_2$ are strictly increasing and that $h_2(x)\to \infty$, as $x \to \infty$. It is now possible to construct a random variable $X$ whose risk function $R$ satisfies
\begin{equation}\label{esimy1}
\limsup_{x \to \infty} \frac{R(x)}{h_1(x)}=1
\end{equation}
and
\begin{equation}\label{esimy2}
\liminf_{x \to \infty} \frac{R(x)}{h_2(x)}=1. 
\end{equation}
First, construct sequence $(x_n)$ by setting $x_0=1$ and $x_{n+1}=h_2^{-1}(h_1(x_n))$ for $n \in \{0,1,2,\ldots\}$.
\begin{figure}
 \centering
 \def\svgwidth{0.8\columnwidth}
 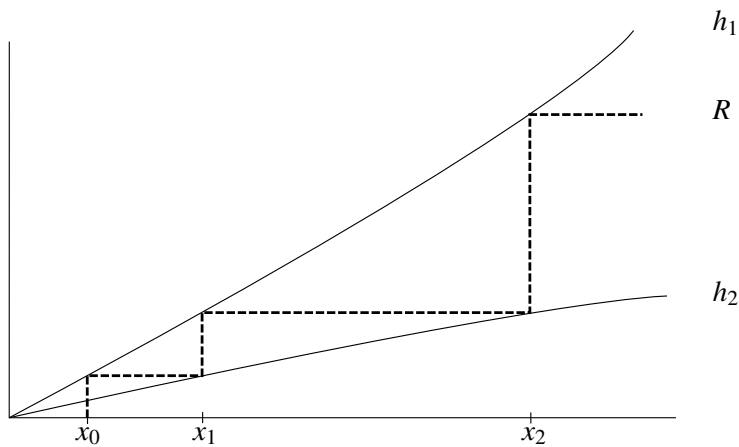
 \caption{Illustration of the hazard function $R$ is drawn using a dashed line. \label{fig1} }
\end{figure}
Now, $x_n \to \infty$, as $n \to \infty$. Moreover, a discrete random variable with a point density function
$$P(X=x_n)= e^{-h_2(x_n)}-e^{-h_1(x_n)}, $$
where $n \in \{0,1,2,\ldots\}$ satisfies requirements \eqref{esimy1} and \eqref{esimy2}. Graphical representation of this construction is given in figure \ref{fig1}.

It is worthwhile to notice that the behaviour of \eqref{esimy2} defines the integrability properties of the random variable $X$, that is, the function $h_2$ is solely responsible of the integrability of $X$. In illustration \ref{fig1} a situation where $h_1$ is concave and $h_2$ is convex is shown. This illustration depicts the fact that a natural scale can be easier to handle than the original function $R=-\log \overline{F}$ itself. 
\end{ex}

Even if the function $R$ is smooth, there may be a better choice of natural scale. This phenomenon can be seen in the following example in the case of log-normal type tails.

\begin{ex}[Two difficult cases in classical theory]\label{esim2}

Two distributions that escape the scope of indices \eqref{e} and \eqref{m} are the Weibull distribution and log-normal type distributions. Weibull distribution is concentrated on $[0,\infty)$ and its tail function has the form 
$$\overline{F}(x)=e^{-\lambda x^{\alpha}}$$
with $\lambda>0$ and $\alpha \in (0,1)$. This is the distribution of $X^{(1/\alpha)}$ when $X$ is exponentially distributed with parameter $\lambda$.  We say that a random variable is of \emph{log-normal type}, if its tail satisfies 
$$\overline{F}(x)\sim cx^\beta e^{-\lambda (\log x)^\gamma},$$
as $x \to \infty$. Here $\beta\in \mathbb{R},\lambda>0$, $\gamma>1$ and $c$ is a positive norming constant.

Let $B$ denote a random variable having Weibull distribution or log-normal type distribution. It is easily seen that $\mathcal{E}(B)=0$ and $\I(B)=\infty$. Hence, the classical indices reveal little information about the distribution. However, since $B$ is heavy-tailed, remark \ref{cor1} ensures that there is a natural scale of the distribution of $B$ that satisfies condition \eqref{summa}. Suppose $(B,B_1,\ldots,B_n)$ are IID variables. Now, the tail of the sum $S_n$ is bounded by the tail of a single variable, that is,
$$P(B_1+\ldots+B_n>x)\leq P(B_1>x)^{1-\epsilon} $$
for all $x$ large enough.

Here, a natural scale can be chosen to be $h(x)=\lambda x^\alpha$ for the Weibull distribution and $h(x)=\lambda (\log x)^\gamma$ for the log-normal type distribution respectively. 
\end{ex}

Next, we calculate up to a constant factor a scale function for a product of two Weibull distributed random variables.

\begin{ex}[Product of two independent Weibull variables]\label{2prod}
Suppose $X$ and $Y$ are IID Weibull distributed random variables with tail $\overline{F}(x)= e^{-\lambda x^{\alpha}}$, where $\lambda>0$ and $\alpha \in (0,1)$. By selecting $g(x_1,x_2)=x_1 x_2$ in Theorem \ref{vtheorem} we get 
$$h_{XY}(x)=c(R_{X}\circ g_d^{-1})(x), $$
where $c>0$ is a constant and $g_d^{-1}(x)=\sqrt{x}$. Since $XY$ has a natural scale of the form $d x^{\alpha/2}$, where $d>0$, we see that the calculation rule
\begin{equation*}
\Ih(XY)=\min(\Ih(X),\Ih(Y))
\end{equation*}
for IID variables $X$ and $Y$ implied by condition \eqref{tulo} of Theorem \ref{laskutheorem} cannot be valid for $h=h_X$. In fact, the scale function of $XY$ grows at a significantly slower speed than the scale function of $X$, which is why the functions $h_{XY}$ and $h_X$ differ by more than a multiplicative constant factor. 
\end{ex}

The previous example dealt with the product of two random variables. The following example shows that Theorem \ref{vtheorem} can be used to obtain simple and general bounds in situations where structure of the model is based on product of IID variables. We study the utility of an economic agent under random IID endowments of commodities using the celebrated Cobb-Douglas model. Different ways of introducing randomness into Cobb-Douglas model with deeper discussion can be found in e.g. \cite{nummelin2000existence}.

\begin{ex}[Tail asymptotics in an economic model based on product structure]\label{tasymp}
Suppose 
$$g(x_1,x_2,\ldots,x_n)=x_1^{a_1}x_2^{a_2}\ldots x_n^{a_n},$$
where $a_1+a_2+\ldots+a_n=1$ and $a_i\geq 0$ for all $i \in \{1,2,\ldots,n\}$. Let $X_1,X_2,\ldots,X_n$ be positive IID variables with common continuous distribution function $F$. Now, $g_d(x)=g_d^{-1}(x)=x$ and, for a given $\epsilon>0$, application of Theorem \ref{vtheorem} yields 
\begin{equation}\label{econarv}
P(g(X_1,X_2,\ldots,X_n)>x)\leq P(X_1>x)^{1-\epsilon}
\end{equation}
for all $x$ large enough.

The function $g$ can be interpreted as a utility function of an economic agent in Cobb-Douglas model. Formula \eqref{econarv} shows that the tail of the utility in random IID allocation of goods is dominated by the tail of a single variable. 
\end{ex}

Next, we move on to the study of the process $(Y_n)$. The following example expands previously knows results to the case of even heavier tails.

\begin{ex}[On asymptotics of tail $\overline{F}_{Y_n}$]\label{esim3}
Suppose $(A_i)$ and $(B_i)$ are independent sequences of positive random variables. In many different fields, such as in insurance and queuing theory, randomly weighted random sums of the type \eqref{disc} appear constantly, see e.g. \cite{nyrhinen2001finite} for background. In this context we set $Y_n=0$.

A variable of interest is 
$$\bar{Y}_n:=\sup_{1\leq k\leq n} Y_k.$$ 
In \cite{Tang2003299} Theorem 4.1. the moment index of the random variable $\bar{Y}_n$ is solved:
$$\I(\bar{Y}_n)=\min(\I(A),\I(B)).$$ 
Using the theorems from previous chapter it is possible to extend this result beyond the scope of polynomial decay. Namely, if the scale function $h$ satisfies properties \eqref{summa} and \eqref{tulo}, a deduction similar to that of \cite{Tang2003299} can be generalised. 

We recall (e.g. from \cite{Tang2003299}) that the process $(\bar{Y}_n)$ admits the recursive representation $\bar{Y}_n\eqdist U_n$, where
$$U_0=0, U_{n+1}=B_{n+1}+A_{n+1}U_n $$
for $n \in \mathbb{N}$ and $\eqdist $ signifies equality in distribution. Here $U_n$ is independent of vector $(A_{n+1},B_{n+1})$. From this it is clear that 
\begin{eqnarray*}
\Ih(\bar{Y}_n)=\Ih(U_n)&=& \Ih(B_{n}+A_{n}U_{n-1}) \\
&=&\min(\Ih(B_n),\Ih(A_n),\Ih(U_{n-1})) \\
& \vdots & \\
&=& \min(\Ih(B),\Ih(A)).
\end{eqnarray*}

Interpretation of formula 
\begin{equation}\label{vikaviit}
\Ih(\bar{Y}_n)=\min(\Ih(B),\Ih(A))
\end{equation}
is that the risk with a more slowly increasing hazard function determines the asymptotic behaviour of the process $(\bar{Y}_n)$. Heuristically, the condition \eqref{tulo} is satisfied when the hazard function grows slower than the logarithmic function, which is not the case with Weibull distributed random variables. This is, however,  the case with log-normal type distributions. Weibull distribution fails to satisfy \eqref{tulo} and example \ref{esim2} shows that property \eqref{vikaviit} cannot be valid for this distribution.
\end{ex}

\begin{ex}[Special case of disconted sum $Y_n$: $B=1$]\label{aesim}

In example \ref{esim3} conditions \eqref{summa} and \eqref{tulo} were taken as assumptions. However, the scale of random variable $Y_n$ is possible to deduce from theorems \ref{vertailuthm} and \ref{vtheorem} alone. 

Consider the discounted sum $Y_n$, where $n\geq 3$. We will make the simplifying assumption $B_i=1$ for all $i \in \{1,2,\ldots,n\}$. In addtition, the sequence $(A_i)$ is assumed to be IID sequence of positive variables with a common continuous distribution function. Using Theorem \ref{vertailuthm} we see that the scale of $Y_n$ is determined by the heaviest of the summands of 
$$Y_n=1+A_1+A_1 A_2+\ldots + A_1\ldots A_{n-1}, $$
which is $A_1\ldots A_{n-1}$. Now Theorem \ref{vtheorem} shows that a natural scale for $Y_n$ is up to a positive constant the function $x \mapsto h_A(x^{1/(n-1)})$. Using Lemma \ref{maxlemma} we see that this scale is also a natural scale for $\bar{Y_n}$.
\end{ex}

The final two examples are of more theoretical nature. Following example shows how a deterministic transformation of a random variable can be used to alter the moment index. For stochastic way to change the moment index the reader is advised to see \cite{oma1}.

\begin{ex}[The scale $h=-\log R$]\label{esim4}

For any random variable $X$ with continuous distribution function, by Theorem \ref{theorem1}, 
$$\I(1/\overline{F}_X(X))=1$$
holds. This way it is possible to give a deterministic transformation, depending on $X$ itself, with which $X$ may be transformed to a random variable with any positive moment index. If the required moment index is $\alpha>0$, the transformation 
$$x \mapsto (1/\overline{F}_X(x))^{(1/\alpha)}$$
can be used. 
\end{ex}

The last example gives a sufficient condition for the moment determinacy of a general non-negative heavy-tailed random variable. For background on the moment problems reader is advised to see \cite{gut2002moment}. The condition is purely a tail condition, that is, modifications of the distribution on a finite interval of $[0,\infty)$ do not change the result.

\begin{ex}[Moment determinacy via decay speed of scale function]\label{esim5} 
Suppose $X\geq 0$ is a heavy-tailed random variable and let $h$ be its natural scale.
We recall that a random variable $X\geq 0$ is determined by its moments if 
\begin{equation*}
E(X^k)=E(Y^k), \, \forall k \in \mathbb{N} \Longrightarrow X\eqdist Y.
\end{equation*}
Using the Hardy's condition it is possible to give a limit test for the moment determinacy involving the concept of natural scale. This test has two benefits compared to other tests such as Carleman condition or finiteness of the Krein integral (see \cite{gut2002moment} for these tests).
\begin{enumerate}
\item Decision is done using the asymptotic properties of the tail function: small realisations of $X$ are irrelevant.
\item No assumption about the absolute continuity w.r.t. Lebesgue measure is needed.
\end{enumerate}

Suppose 
\begin{equation}\label{momdet}
\liminf_{x \to \infty} \frac{h(x)}{\sqrt{x}}>0.
\end{equation}
Then the distribution of $X$ is determined by its moments. Implication is easily verified by observing that \eqref{momdet} together with the definition of natural scale and positivity yields
$$\liminf_{x \to \infty} \frac{R(x)}{\sqrt{x}}\geq \left( \liminf_{x \to \infty} \frac{h(x)}{\sqrt{x}} \right) \left(\liminf_{x \to \infty} \frac{R(x)}{h(x)}\right) >0.$$
Hence, by connection of Theorem \ref{theorem1}, there exists $c>0$ such that 
$$E(e^{c\sqrt{X}})<\infty $$
and Theorem 1 of of \cite{stoyanov2012hardy} confirms that $X$ is determined by its moments.
\end{ex}

\begin{remark} In example \ref{esim5} $X$ is assumed heavy-tailed. This is not a limitation, because any light-tailed random variable would automatically satisfy condition 
$$E(e^{cX})<\infty$$
for some $c>0$ and instantly be determined by its moments.
\end{remark}

\appendix

\section{Omitted technical details}\label{mappendix}

\subsection{Construction of the function $h_\eta$ of Theorem \ref{theorem1}}\label{eappendix}

In Theorem \ref{theorem1} it was claimed that a function satisfying \eqref{tuplayht} exists. To construct function $h_\eta$, recall that the function $h$ itself is assumed increasing and continuous. Therefore, we may construct sequence $(y_k)$, where $y_k$ is defined to be the unique solution of equation $h(x)=k \eta$. 

Now, we set $h_\eta(y_k)=k \eta$ for all $k\in \mathbb{N}$ and define families of functions $\mathcal{A}_k$ in the following way: $f \in \mathcal{A}_k$ if and only if 
\begin{enumerate}
\item $f\colon [y_k,y_{k+1}]\to \mathbb{R}$ is a concave function
\item $f(y_k)=k\eta$ and $f(y_{k+1})=(k+1)\eta$ and
\item $f(x)\geq h(x)\, \forall x \in [y_k,y_{k+1}]$.
\end{enumerate}
Each set $\mathcal{A}_k$ is not empty since $h$ is continuous. Now, if $x \in (y_k,y_{k+1})$ we set
$$h_\eta(x)=\inf_{f \in \mathcal{A}_k}\{f(x)\}.  $$
By construction \eqref{tuplayht} is valid.

\subsection{Construction of the function $g$ of Theorem \ref{theorem3}} \label{tappendix}

Suppose we are given a hazard function $R$ of a heavy-tailed random variable. We must now find a continuous function $g$, depending on $R$, such that $g(x)=o(R(x))$, as $x \to \infty$.

Since $R$ is right continuous, we may define sequence $(x_n)$ by $x_0:=0$ and 
\begin{equation*}
x_n:=\min \{x:R(x)\geq n\}, \, n \in \mathbb{N}.
\end{equation*}
We note that $x_n\to \infty$, as $n \to \infty$. Define $g(x_0)=g(0)=0$ and
\begin{equation*}
g(x_n):= \sqrt{n-1}
\end{equation*}
for all $n \in \mathbb{N}$. Between points of the sequence $(x_n)$ values of $g$ are given by linear interpolation
\begin{equation*}
g(x_{n+1}-(1-t)(x_{n+1}-x_n))=g(x_{n+1})-(1-t)(g(x_{n+1})-g(x_n)), \, t \in [0,1]
\end{equation*}
for all $n\in \{0,1,2,\ldots\}$.
The function $g$ is continuous and by construction $g(x)=o(R(x))$, as $x \to \infty$.

\subsection{Details for remark \ref{remark1}}\label{rappendix}

Suppose we are given a concave function function of remark \ref{remark1}. By concavity, for any $y\in (0,1)$ and $x>0$,
\begin{equation}\label{convehto}
h(yx)=h(yx+(1-y)0)\geq yh(x)+(1-y)h(0)=yh(x). 
\end{equation}
Using this we get subadditivity: for all $a,b>0$
\begin{eqnarray*}
h(a+b)&=&\frac{a}{a+b}h(a+b)+\frac{b}{a+b}h(a+b)\\
&\leq &h(a)+h(b).
\end{eqnarray*}
Formula \eqref{convehto} implies that $h(yx)/h(x)\leq y$, which is why 
$$\limsup_{x \to \infty} \frac{h(yx)}{h(x)}<\infty $$
holds and $h\in \mathcal{D}$.

\section*{Acknowledgements}

The deepest gratitude is expressed to the Finnish Doctoral Programme in Stochastics and Statistics (FDPSS) and the Centre of Excellence in Computational Inference (COIN) for financial support (Academy of Finland grant number 251170). Special thanks are due to Ph.D. Harri Nyrhinen for numerous enlightening conversations and guidance.

\bibliographystyle{plain}
\bibliography{mybib}

\end{document}